\documentclass[11pt]{amsart}
\usepackage{amscd,amssymb,palatino}
\usepackage[utf8]{inputenc}
\usepackage[english]{babel}
\usepackage{caption}
\usepackage{etex}
\usepackage[leqno]{amsmath}
\usepackage{amssymb}
\usepackage{color}
\usepackage{shadow}
\usepackage{epsfig}
\usepackage{epic}
\usepackage{eepic}
\usepackage{graphics}
\usepackage{graphicx}
\usepackage{psfrag}
\usepackage{calc}

\usepackage{tikz}
\usepackage[dvipsnames,prologue,table]{pstricks}
\usepackage{pstcol}
\usetikzlibrary{arrows,decorations,patterns,positioning,automata,shadows,fit,shapes,calc,decorations.markings,backgrounds,scopes,decorations.text}

\usepackage{tipa}

\usepackage{fancyhdr}
\fancypagestyle{plain}{\fancyhf{}
\fancyfoot[C]{{\large\sf\bfseries\thepage}}}
 \DeclareMathOperator{\rank}{rank}

\usepackage{calc}
\newcommand\1{\lower 9pt\hbox{\underbar{}}}
\numberwithin{equation}{section}

\newtheorem {Theorem}                   {Theorem}

\newtheorem {RefTheorem}[equation]      {Theorem}

\theoremstyle{definition}
\newtheorem {Definition}[equation]{Definition}
\newtheorem {Remark}[equation]          {Remark}

\newcommand{\pr} {\smallskip\noindent{\bf Proof\,\,}}

\begin{document}

\title[Quantization of $b^m$-symplectic manifolds]{On geometric quantization of $b^m$-symplectic manifolds}

\author{Victor W. Guillemin}\thanks{V. Guillemin is supported in part by a Simons collaboration grant.}
\author{Eva Miranda}
\thanks{{ E. Miranda  is supported by the Catalan Institution for Research and Advanced Studies via an ICREA Academia Prize 2016, a Chaire d'Excellence de la Fondation
Sciences Math\'{e}matiques de Paris and partially supported  by the grants reference number MTM2015-69135-P (MINECO/FEDER) and reference number 2014SGR634 (AGAUR). This work
is supported by a public grant overseen by the French National Research Agency (ANR) as part of the \emph{\lq\lq Investissements d'Avenir"} program (reference:
ANR-10-LABX-0098).}}
\author{Jonathan Weitsman}
\thanks{J. Weitsman was supported in part by NSF grant DMS 12/11819 and by a Simons Collaboration Grant}
\address{Department of Mathematics, MIT, Cambridge, MA 02139}
\email {vwg@math.mit.edu}
\address{{Department of Mathematics}, Universitat Polit\`{e}cnica de Catalunya and BGSMath, Barcelona, Spain \\ CEREMADE (Universit\'{e} de Paris Dauphine), IMCCE
(Observatoire de Paris) and IMJ (Unversit\'{e} de Paris 7), 77, avenue Denfert Rochereau
75014, Paris, France}
\email{eva.miranda@upc.edu, Eva.Miranda@obspm.fr}
\address{Department of Mathematics, Northeastern University, Boston, MA 02115}
\email{j.weitsman@neu.edu}
\thanks{\today}

\begin{abstract} We study the formal geometric quantization of  $b^m$-symplectic manifolds equipped with Hamiltonian actions of a torus $T$ with nonzero leading modular
weight.  The resulting virtual $T-$modules are finite dimensional when $m$ is odd, as in \cite{gmw2}; when $m$ is even, these virtual modules are not finite dimensional, and
we compute the asymptotics of the representations for large weight.


\end{abstract}
\maketitle

\section{Introduction}

The purpose of this paper is to continue our study of formal geometric quantization for a class of Poisson manifolds, the $b^m$-Poisson manifolds which appeared first in the
thesis of Geoffrey Scott \cite{scott}, in the case where these manifolds are equipped with Hamiltonian torus actions. In the case $m=1$, those are  the $b$-symplectic
manifolds of \cite{guimipi, guimipi2} whose quantizations are constructed in \cite{gmw2}.  The methods we use here are very similar to those of \cite{gmw2}, and as in the
case of $b$-symplectic manifolds, the quantizations turn out to have remarkable properties:  We obtain finite dimensional virtual modules for $m$ odd and infinite dimensional
virtual modules with remarkable asymptotic properties for $m$ even.  We hope this repertory of examples will give further intuition about what might be a geometric
quantization for more general Poisson manifolds.

Let $m$ be a positive integer. A $b^m$-symplectic manifold is a smooth manifold $M$, along with a smooth hypersurface {$Z\subset M$}, and a choice of an $m$-germ of a
$\mathcal{C}^{\infty}$-function at $Z$ along with a closed, nondegenerate $b^m$-form  $\omega$ of degree $2$ on $M$ (see \cite{scott} for details; we give a brief summary of
Scott's results below).

As in \cite{gmw2}, as a form of prequantization, we assume that the $b^m$-symplectic form is {\em integral,} in the sense that this form has integral decomposition under
Scott's analog of the {Mazzeo--Melrose} theorem.

Suppose that $(M,\omega)$ is compact, oriented, and connected and is equipped with a Hamiltonian action of a compact torus $T$, with nonzero leading modular weight (see \cite{gmw3}).  Suppose we are given a line bundle $L \to M$ and a connection $\nabla$ on $L|_{M-Z}$ whose curvature is $\omega|_{M-Z}.$\footnote{Such a line bundle with connection is needed for the formal quantization to be well-defined, even in the $b$-symplectic case discussed in \cite{gmw2}.  We would like to thank the referee for pointing out this omission in \cite{gmw2} and in an earlier version of this paper.}

We then define the formal geometric quantization $Q(M)$\footnote{We write $Q(M)$ by abuse of notation even though the quantization as defined depends on the data $\omega, L, \nabla.$} as a virtual $T$-module by a procedure analogous to that in \cite{Wei, gmw2}; then for each weight $\alpha \in
\mathfrak{t}^*$ of the torus $T,$ the multiplicity of the weight $\alpha$ in $Q(M)$ is given by

$$Q(M)^\alpha=\pm Q(M//_{\alpha} T)$$

\noindent where $Q(M//_{\alpha} T)$  is the geometric quantization of the symplectic manifold $M//_{\alpha}T$\footnote{Some care must be taken when $\alpha$ is a singular
value of the moment map; see below.} and the sign is positive if the symplectic orientation on the symplectic quotient $M//_{\alpha}T$ agrees with the orientation inherited from $M$, and negative otherwise.\footnote{This sign convention is inspired by the results of \cite{ckt} for presymplectic manifolds; in the symplectic case, these orientations always agree, and the sign is always positive.}

Our main theorem is the following:

\begin{Theorem}\label{main} Suppose $T$ acts on $M$ with nonzero leading modular weight. Then,

\begin{enumerate}

\item If $m$ is odd, $Q(M)$ is a finite dimensional virtual $T$-module.

\item If $m$ is even, there exists a weight $\xi\in\mathfrak{t}^*$, integers $c_{\pm}$, and $\lambda_0>0$ such that if $\lambda>\lambda_0$, and $\eta\in\mathfrak{t}^*$ is a
    weight of $T$,

    $$\dim{Q(M)^{\lambda\eta}}=\begin{cases} 0 \quad if \quad  \eta\neq \pm\xi \\ c_{\pm} \quad if \quad  \eta= \pm\xi\end{cases}$$

    (In fact, $c_{\pm}= \epsilon_\pm \dim{Q(M)^{\pm\lambda\xi}}, $ where $\epsilon_\pm \in \{ \pm 1 \},$ for any $\lambda$ sufficiently large.)

    \end{enumerate}
\end{Theorem}
    \begin{Remark}
    In the case of $m=1$, Part (1) of this theorem was proved in \cite{gmw2}. As in that case, it is natural to conjecture there exists a natural Fredholm operator  with
    index given by $Q(M)$.
    \end{Remark}

    \begin{Remark}
    Unlike in the case of $m=1$, however, there is no reason to expect the signed symplectic volume of $(M, \omega)$ to be finite. Thus the naive version of the semiclassical
    limit may not hold.

    \end{Remark}

{\textbf{Acknowledgements:}  We thank C\'{e}dric Oms for carefully reading a first version of this article. We would also like to thank the referee for several helpful comments and corrections.  In particular, the referee pointed out to us the delicate issue of the potential dependence of quantization on the choice of a connection with curvature given by the symplectic form, in the case where the quotients may be orbifolds.  This issue was overlooked in \cite{gmw2}, and we address it in this paper.  It does not affect the statements of the main theorems and addressing it requires only some care with the proofs, not any major changes.

\section{$b^m$-Manifolds}

Let $M$ be a compact manifold, and let $f\in\mathcal{C}^\infty(M)$ have a transverse zero at  a hypersurface $Z\subset M$. Let $m$ be a positive integer, the $m$-germ of $f$
at $Z$ gives rise to a locally free sheaf, and therefore a vector bundle, given by,

$$\Gamma(^{b^m}TM)=\{v \in \Gamma(TM): vf \quad \text{vanishes to order $m$ at $Z$} \}.$$

\noindent  By considering sections of the wedge powers ${\Lambda^k} (~^{b^m}{T^*M})$ we obtain a complex $(~^{b^m}\Omega^k(M), d)$ of differential forms with singularities at
$Z$ and an associated de Rham cohomology (the so-called $b^m$-cohomology).

\begin{RefTheorem}[\textbf{$b^m$-Mazzeo--Melrose}, \cite{scott}]\label{thm:Mazzeo-Melrose}
\begin{equation}\label{eqn:Mazzeo-Melrosebm}
  ~^{b^m}H^p(M) \cong H^p(M)\oplus(H^{p-1}(Z))^m.
\end{equation}
\end{RefTheorem}

We call an element of $~^{b^m}H^p(M)$ {\em integral} if its image under the isomorphism above is integral.

A two form $\omega\in~^{b^m}\Omega^2(M)$ is $b^m$-symplectic if it is closed and nondegenerate (as an element of ${\Lambda^2}(~^{b^m}T^*M$)). We assume now that $M$ is
equipped with an integral $b^m$-symplectic form.  This is an analog for $b^m$-symplectic manifolds of the prequantization condition of \cite{gmw2}.

In this case, the image of $\omega$ in $H^2(M)$ gives rise to a line bundle on $M$. In order to study the meaning of integrality of the other terms in the Mazzeo--Melrose
formula (\ref{eqn:Mazzeo-Melrosebm}), we recall that in a neighborhood $U$ of the critical set $Z$, $U=Z\times (-\epsilon, \epsilon)$,   $\omega$ may be written as

\begin{equation}\label{eqn:newlaurent}
{\omega = \sum_{j = 1}^{m}\frac{df}{f^j} \wedge \pi^*(\alpha_{j}) +  \beta }
\end{equation}
\noindent where the $\alpha_j$ are closed one forms on $Z$, $\beta$  is a  closed 2-form on $U$, and $\pi:U\longrightarrow Z$ is the projection. Nondegeneracy of $\omega$
implies that $\beta\vert_{Z}$ is of maximal rank and $\alpha_m$ is nowhere vanishing. In fact,   $\alpha_m$ defines the symplectic foliation of the Poisson structure
associated to $\omega$, and $\beta$ gives the symplectic form on the leaves  of this foliation. In our case, the integrality of $\alpha_m$ gives a map $\phi:Z\longrightarrow
S^1$ whose fibers define the symplectic foliation showing that this foliation has compact leaves and is a mapping torus.

Now let us assume a torus $T^d$ acts effectively on $M$ preserving $Z,f$ and $\omega$, and furthermore that this action is Hamiltonian in the sense of \cite{gmw3}; that is,
there exists a moment map $\mu\in ~^{b^m}\mathcal{C}^\infty(M)\otimes \mathfrak{t}$ with

$$\langle d\mu, \xi \rangle= i_{\xi^M}\omega$$

\noindent for any $\xi\in\mathfrak{t}$ and where $\xi^M$ stands for the fundamental vector field generated by $\xi$;

\noindent where $$~^{b^m} \mathcal{C}^\infty(M)= f^{-(m-1)}\mathcal{C}^\infty(M) + ~^{b} \mathcal{C}^\infty(M)$$

\noindent and $$^{b}\mathcal{C}^\infty(M)=\{ g \log\vert f\vert+ h, g, h\in \mathcal{C}^\infty(M)\}.$$

As in \cite{gmw3}, the action gives rise to a series of modular weights $a_1, \dots, a_m\in \mathfrak{t}^*$ in each connected component of $Z$ given by

$$a_j(\xi)=\alpha_j(\xi^M);$$

\noindent in \cite{gmw3} we show these are constants.

We now assume

\textbf{Assumption:} $a_m$ is nonzero.

As in \cite{gmw3} it is sufficient to assume $a_m$ is nonzero for one connected component of $Z$.

Let $\mathfrak{t}_L = (a_m)^0 \subset \mathfrak{t}$ be the annihilator of the weight $a_m$ of $T$ and let $T_L$ be the torus generated by $\mathfrak{t}_L.$  This is a one
codimensional subtorus of $T.$


Then the main technical result of \cite{gmw3} is the following local convexity theorem:

\begin{RefTheorem}[\textbf{Guillemin-Miranda-Weitsman}\label{thm:localconvexity}]
There exists a neighborhood $U=Z\times (-\epsilon, \epsilon)$ where the moment map $\mu: M - Z\longrightarrow\mathfrak{t}^*$ is given by

{$$\mu=\sum_{i=1}^{m-1}\frac{1}{g^i} a_{i+ 1} +\log \vert g \vert a_1+\mu_L$$}

\noindent where $g$ is a function agreeing with $f$  to order $m$, {$a_i\in\mathfrak{t}^0_L$}, and $\mu_L$ is the moment map for the $T_L$-action on the symplectic leaves of
the foliation.

\end{RefTheorem}

\section{Formal Geometric Quantization}

We wish to define geometric quantization for $b^m$-symplectic manifolds, by analogy with some form of geometric quantization in the symplectic case.  Ideally this would involve constructing a Fredholm operator on the manifold, computing its index, and proving a vanishing theorem, so that we would obtain a vector space.  Since we do not know how to perform this construction, we try to obtain intuition about what such a quantization would look like by using the expectation that quantization commutes with reduction to define a "formal geometric quantization" in the $b^m$-symplectic case, in analogy with a similar construction \cite{Wei,p} in the case of noncompact Hamiltonian $T$-spaces with proper moment map.

Let $(M,\omega)$ be an integral $b^m$-symplectic manifold
with a Hamiltonian action of a torus $T:= \mathbb T^d$ with nonzero leading modular weight.  Suppose we are given a line bundle $L\to M$ and a connection $\nabla$ on $L|_{M-Z}$ with curvature $\omega|_{M-Z}.$

\begin{Definition}[\cite{Wei}]\label{deffgq}  A virtual $T$-module $V$ is the {\em quantization} of $M$ (and we write $Q(M)=V$) if, for every compact integral Hamiltonian
$T$-space $N$ we have

\begin{equation}\label{qreqn1}
(V\otimes Q(N))^T = \epsilon Q((M \times N)//_0T)\end{equation}

\noindent where $Q(N)$ denotes the standard geometric quantization of $N$, and $Q(M \times N)//_0T$ is the geometric quantization of the compact integral symplectic manifold
$(M \times N)//_0T,$, i.e., the symplectic reduction of $M \times N$ at $0$\footnote{Note that since the moment map is singular on $Z,$ and $N$ is compact, $(M \times N)//_0
T = ((M-Z) \times N)//_0 T $, so that $(M \times N)//_0 T $ is compact and symplectic.}, and where $\epsilon$ is $+1$ if the symplectic orientation on the symplectic quotient $(M \times N)//_0T$ agrees with the orientation inherited from $M \times N$, and $-1$ otherwise.\footnote{This sign convention is inspired by the results of \cite{ckt} for presymplectic manifolds; in the symplectic case, these orientations always agree, and the sign is always positive.}

\end{Definition}

{See also \cite{p, vergne1}}.
\begin{Remark}
In the case where $(M \times N)//_0T$ is not smooth, the quantization must be defined by the shift desingularization of Meinrenken and Sjamaar \cite{ms}, which is standard by
now (see e.g. \cite{Wei,p}).

\end{Remark}

Equation (\ref{qreqn1}) means that $Q(M)=\oplus_i \epsilon_i Q((M \setminus Z)_i),$ where the $(M\setminus Z)_i$ are the components of $M \setminus  Z$, where $Q(M \setminus Z)$ is the formal geometric quantization of the non-compact Hamiltonian $T$-space (with proper moment
map) $M \setminus Z$ whose quantization was defined in \cite{Wei,p}, and where $\epsilon_i \in \{\pm 1\}$ are determined by the relative orientations of the symplectic forms on the components of $M \setminus Z$ and the overall orientation of $M.$  (Integrality of the symplectic manifold $M \setminus Z$ follows from our integrality condition on the $b^m$-symplectic form, and properness of the moment map follows from the condition that the leading modular weight is nonzero.)  Thus the quantization of $M$ exists and is unique.  Alternatively,
$$Q(M) = \bigoplus_\alpha \epsilon(\alpha) Q(M//_\alpha T) \alpha,$$

\noindent where $Q(M//_\alpha T)$ must again be defined using the shifting trick where $\alpha$ is not a regular value of the moment map, and $\epsilon(\alpha) \in \{ \pm 1\}$ is determined by the relative orientations of $M$ and $M//_{\alpha} T.$

\section{Proof of the main theorem}

Since $$Q(M) = \bigoplus_\alpha \epsilon(\alpha) Q(M//_\alpha T) \alpha,$$ \noindent
it suffices to compute $Q(W)$ where $W$ is some $T$-invariant neighborhood of $Z$.
 The complement
of this neighborhood is compact, and therefore contributes only a finite number of weights to $Q(M)$; so it does not affect either the finiteness statement (1) or the asymptotic statement (2).

{Consider, then, a component $Z_i$ of $Z$ and a neighborhood $U=(-\epsilon, \epsilon)\times Z_i$ containing $Z_i$.}  The local convexity theorem of \cite{gmw3} (see Theorem \ref {thm:localconvexity})  shows that the moment map in the region $U$ may be taken as

{$$\mu(x,z)=\sum_{i=1}^{m-1}\frac{a_{i+1}}{x^i}+ a_1 log\vert x\vert + \mu_L(p(z))$$}

\noindent where $p:Z\longrightarrow L$ is the projection onto the fiber of the trivial symplectic fibration $\phi: Z\longrightarrow S^1$.

 In order to prove the second part of the theorem, observe that any weight can be written as $\xi=\xi_0+\eta\in \mathfrak{t}^*$ where $\xi_0$ lies in the one dimensional
 subspace of $\mathfrak{t}$ generated by $a_m$ and $\eta\in  \mathfrak{t}_L^*$.

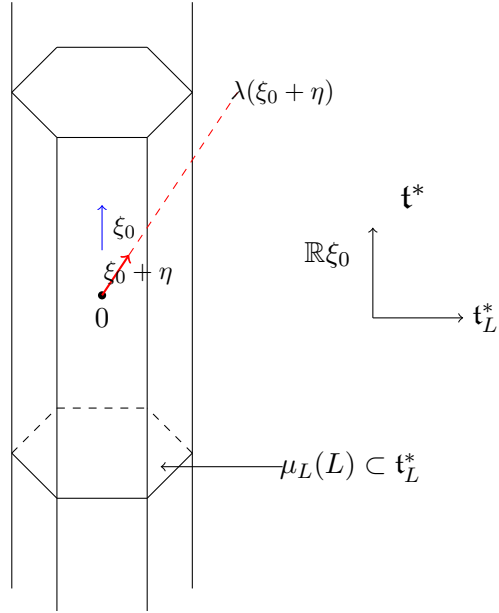
 \begin{center}
\begin{tikzpicture}[scale=0.6]

\draw[-] (-1,0,0)--(1,0,0);
\draw[-] (1,0,0) -- (2,1,0);
\draw[-] (-1,0,0) -- (-2,1,0);
\draw[-] (-2,-2,0) -- (-2,11,0);
\draw[-] (2,-2,0) -- (2,11,0);
\draw[dashed] (-2,1,0) -- (-1,2,0);
\draw[dashed] (2,1,0) -- (1,2,0);
\draw[dashed] (-1,2,0) -- (1,2,0);

\draw[-] (2,9,0)-- (1,8,0);
\draw[-] (-2,9,0)-- (-1,8,0);
\draw[-] (-1,8,0)-- (1,8,0);
\draw[-] (2,9,0)-- (1,10,0);
\draw[-] (1,10,0)-- (-1,10,0);
\draw[-] (-2,9,0)-- (-1,10,0);

\draw[-] (-1,8,0) -- (-1,-2.5,0);
\draw[-] (1,8,0) -- (1,-2.5,0);

\draw[->] (4,0.7,0) -- (1.3,0.7,0);

\draw (5.5,0.7,0) node [] {$\mu_L(L) \subset \mathfrak{t}_L^*$};

\draw[blue][->] (0,5.5,0) -- (0,6.5,0);
\draw (0.5,6,0) node[scale=0.9] {$ \xi_0 $};
\draw (0,4.5,0) node[circle,fill, scale=0.3]{};
\draw (0,4,0) node[]{$0$};
\draw[red][->, thick] (0,4.5,0)-- (0.6,5.4,0);
\draw (0.8,5,0) node [scale=0.9] {$\xi_0 + \eta$};
\draw[red][dashed] (0.6,5.4,0) -- (3,9,0);
\draw (4,9,0) node [scale=0.9] {$\lambda (\xi_0 + \eta)$};

\draw[->] (6,4,0) -- (8,4,0);
\draw[->] (6,4,0) -- (6,6,0);
\draw (8.5,4,0) node[] {$ \mathfrak{t}_L^* $};
\draw (5,5.4,0) node[] {$\mathbb{R}\xi_0$};
\draw (6.9,6.7,0) node[scale=1.2] {$ \mathfrak{t}^* $};
\end{tikzpicture}
\end{center}

First we note that $Q(M)^{\zeta}=0$ unless $\zeta \in Im(\mu_L),$ since the reduced space will be empty if $\zeta$ does not lie in the image of $\mu$.
However, for sufficiently large $|\lambda|,$ $\lambda\xi$ is not in the image of the moment map when $\eta\neq0$, thus proving that $Q(M)^{\lambda\xi}=0$ unless $\eta=0$ for
sufficiently large $|\lambda|.$ This  proves the vanishing statement in (2).

To compute the asymptotics of the quantization, we need the following result, in the case where $(M,\omega)$ is a symplectic manifold.

 \begin{RefTheorem}\label{quantinstagesham}  Let $(M,\omega)$ be a (possibly noncompact) Hamiltonian $T-$space equipped with a $T-$ equivariant line bundle $L$ with connection $\nabla$ of curvature $\omega.$  Suppose $T' \subset T$ is a subtorus, and that the moment maps for the $T$ and $T'$ actions on $M$ are both proper.  Write $Q'(M)$ for the formal geometric quantization of $M$ as a Hamiltonian $T'$ space with proper moment map, and note that $Q'(M)$ is a virtual $T$-module.  Then as virtual $T$-modules,

$$Q'(M) = Q(M).$$
\end{RefTheorem}

The proof is an immediate consequence of the fact that quantization commutes with reduction, and has the following corollary:

\begin{RefTheorem}\label{quantinstages}

 Let $(M,\omega)$ be a $b^m$-symplectic manifold equipped with a Hamiltonian $T$ action, and with a $T-$ equivariant line bundle $L$ with connection $\nabla$ on $M \setminus Z$ of curvature $\omega \vert_{M \setminus Z}.$  Suppose $T' \subset T$ is a subtorus, and that both $T$ and $T'$ act on $M$ with nonzero leading modular weight.  Write $Q'(M)$ for the formal geometric quantization of $M$ as a $b^m$-Hamiltonian $T'$ space, and note that $Q'(M)$ is a virtual $T$-module.  Then as virtual $T$-modules,

$$Q'(M) = Q(M).$$

\end{RefTheorem}

We now continue the proof of Theorem \ref{main}. To study the asymptotic behavior of the quantization $Q(M),$ we note that Theorem \ref{quantinstages} means that it suffices to consider the case where $T$ is a circle whose moment map is given by the nonzero modular weight.  In a neighborhood of $Z,$ we can assume without loss of generality that this circle acts freely (by taking a quotient by a finite stabilizer group, if necessary).  This means that the symplectic quotients are manifolds, so that the dimension or their quantizations depend only on the symplectic form, and not on the choice of connection (this follows from the Riemann Roch formula for the dimension of the quantization of a symplectic manifold).

Next, we note that the normal form (\ref{eqn:newlaurent}) for the symplectic form shows that in a neighborhood of a
connected component of $Z,$ the reduced spaces $M//_{\lambda \xi_0} T$ are symplectomorphic to the leaves of the symplectic foliation on $Z. $  In particular, in such a neighborhood, i.e., for sufficiently large $\lambda,$ these symplectic quotients are independent of $\lambda,$ and the dimension of the quantization is given by
$c_{\pm}= \pm \rank{Q(M)^{\pm\lambda\xi}}$ for any $\lambda$ sufficiently large.

 To prove part (1), for $m$ odd,  we note that the explicit formula for the symplectic form in a neighborhood of a component of $Z$ shows that in such a neighborhood, each reduced space consists of two components, one whose symplectic orientation agrees with that inherited from $M,$ and one where the orientations are opposite; and that these reduced spaces are symplectomorphic.  The quantization of such a neighborhood is the trivial $T-$module.

\end{document}